\documentclass{amsart}
\usepackage{amssymb}
\usepackage{amsmath}
\usepackage{amscd}
\usepackage{url}
\usepackage{bbm}
\usepackage{tikz-cd}

\newtheorem{prop}{Proposition}[section]
\newtheorem{prop*}{Proposition*}[section]

\newtheorem{lem*}[prop]{Lemma*}
\newtheorem{thm}[prop]{Theorem}
\newtheorem{thm*}[prop]{Theorem*}

\theoremstyle{definition}

\newtheorem{remar}[prop]{Remark}

\newcommand{\ord}{{\mathrm {ord}}}
\newcommand{\End}{{\mathrm {End}}}
\newcommand{\Hom}{{\mathrm {Hom}}}

\newcommand{\Ind}{{\mathrm {Ind}}}

\newcommand{\Nm}{{\mathrm {Nm}}}

\newcommand{\Tr}{{\mathrm {Tr}}}

\newcommand{\Frob}{{\mathrm {Frob}}}
\def\rank{\mathop{\mathrm{ rank}}\nolimits}
\newcommand{\Gal}{\mathrm {Gal}}

\newcommand{\Res}{{\mathrm {Res}}}

\newcommand{\A}{{\mathbb A}}
\newcommand{\CC}{{\mathbb C}}
\newcommand{\C}{{\mathbb C}}

\newcommand{\R}{{\mathbb R}}

\newcommand{\QQ}{{\mathbb Q}}
\newcommand{\Q}{{\mathbb Q}}
\newcommand{\ZZ}{{\mathbb Z}}
\newcommand{\Z}{{\mathbb Z}}

\newcommand{\TT}{{\mathbb T}}

\newcommand{\HH}{{\mathfrak H}}

\newcommand{\OOO}{{\mathcal O}}

\newcommand{\pp}{{\mathfrak p}}

\newcommand{\n}{{\mathfrak n}}

\newcommand{\q}{{\mathfrak q}}
\newcommand{\p}{{\mathfrak p}}

\newcommand{\mm}{{\mathfrak m}}

\newcommand{\FF}{{\mathbb F}}

\newcommand{\GL}{\mathrm {GL}}

\newcommand{\Qbar}{\overline{\mathbb Q}}

\DeclareFontFamily{U}{wncy}{}
    \DeclareFontShape{U}{wncy}{m}{n}{<->wncyr10}{}
    \DeclareSymbolFont{mcy}{U}{wncy}{m}{n}
    \DeclareMathSymbol{\Sh}{\mathord}{mcy}{"58} 

\begin{document}
\title{Modular elliptic curves and hyperbolic uniformization}
\author{Neil Dummigan}
\address{University of Sheffield\\ School of Mathematical and Physical Sciences\\
Hicks Building\\ Hounsfield Road\\ Sheffield, S3 7RH\\
U.K.}
\email{n.p.dummigan@shef.ac.uk}
\author{Devendra Tiwari}
\email{devendra9.dev@gmail.com}


\date{ August 25th, 2026.}

\begin{abstract} In an article published a few years before the modularity of elliptic curves over $\Q$ was proved, Mazur \cite{maz} looked at modularity as a purely complex analytic phenomenon, defining a notion of an elliptic curve over $\Q$ having a hyperbolic uniformisation of arithmetic type. Such an elliptic curve (of conductor $N$, say) is necessarily geometrically modular, i.e. a quotient of the jacobian of the modular curve $X_0(N)$, by a morphism defined over $\Q$. We extend these ideas to elliptic curves over totally real fields of odd degree, using Shimura curves for quaternion algebras split at all finite places and one real place. In particular, we prove that the existence of a hyperbolic uniformisation of arithmetic type would imply geometric modularity.
\end{abstract}

\maketitle

\section{Introduction} 
We begin by introducing enough notation and background to give a precise statement of the main theorem. 
Let $F$ be a totally real number field, of {\em odd} degree $[F:\Q]=d$. Let $\infty_1,\ldots,\infty_d$ be the real places of $F$, with associated embeddings $\tau_i:F\rightarrow\R$, and let $B/F$ be the quaternion algebra that is split at $\infty_1$, ramified at $\infty_2,\ldots,\infty_d$ (an even number of places) and split at all finite places. Thus, letting $F_v$ be the completion of $F$ at a place $v$, and $B_v:=B\otimes_F F_v$, $B_v$ is a division algebra for $v=\infty_2,\ldots,\infty_d$, while we may fix isomorphisms $\iota_v: B_v\simeq M_2(F_v)$ for all other places.  Let $\nu: B^{\times}\rightarrow F^{\times}$ be the reduced norm, and use the same symbol to denote its extensions to all the $B_v^{\times}$, to $B^{\times}_{\A}:=B^{\times}(\A)$, where $\A=\A_F$ denotes the ring of adeles of $F$, and to its infinite and finite components $B^{\times}_{\infty}$ and $B^{\times}_f$ respectively. Let the index-$2$ subgroups $B^{\times}_+$ and $B^{\times}_{\infty +}$ be defined by the condition $\tau_1(\nu(b))>0$, and $F^{\times}_+$ by $\tau_i(x)>0\,\,\,\forall\, 1\leq i\leq d$. (Note that automatically $\tau_i(\nu(b))>0$ for $2\leq i\leq d$.) For a comprehensive treatment of the local and global properties of quaternion algebras over number fields, and their Eichler orders, we refer the reader to \cite[Ch 14, 27 and 28]{jv}.

Let $U$ be an open compact subgroup of $B^{\times}_f$. There is an associated Shimura curve $X_U/F$, cf. \cite[(1.1.1),(1.1.2)]{car}. This is a proper, smooth curve defined over $F$, connected but not in general geometrically connected. Let $S=S(U):=B^{\times}_{\infty +}U\leq B^{\times}_{\A +}$. Let $F_U$ be the abelian extension of $F$ corresponding via class field theory to the open subgroup $F^{\times}\nu(S)=F^{\times}F^{\times}_{\infty +}\nu(U)$ of $\A^{\times}$. By global class field theory, the open subgroup $F^\times \nu(S)$ of the idèle group $\mathbb{A}_F^\times$
uniquely determines a finite abelian extension $F_U / F$ via the Artin reciprocity map (see \cite[Ch. VI]{ne} for more details). If $B^{\times}_{\A +}$ is a disjoint union $\cup_{i=1}^h B^{\times}_+ x_iS$ (let's say $x_1=1$) and $S_i:=x_iSx_i^{-1}$ then the geometric components of $X_U$ are defined over $F_U$, and are $V_{S_i}$ in Shimura's notation \cite[(2), 2.5]{shi70}, \cite[Theorem 2.3.1]{oh}. They are conjugate under the action of $\Gal(F_U/F)$, with the Artin symbol $[\nu(x_i^{-1}), F]$ taking $V_{S}$ to $V_{S_i}$. The Shimura curves associated with $B_f^{\times}$-conjugate open compact subgroups are isomorphic over $F$.

Each component $V_{S_i}$, as a geometrically connected curve over $F_U$, has a jacobian $\mathrm{Jac}(V_{S_i})$. This is an abelian variety defined over $F_U$ such that, for any extension $k/F_U$ in which $V_{S_i}$ has a point, $\mathrm{Jac}(V_{S_i})(k)\simeq \mathrm{Pic}^0((V_{S_i})_k)$ (the group of degree-zero $\Gal(\Qbar/k)$-invariant divisor classes). The Shimura curve $X_U$ has something we call its jacobian variety $J_U$, an abelian variety defined over $F$. It may be obtained from any of the $\mathrm{Jac}(V_{S_i})$ by Weil restriction of scalars from $F_U$ to $F$. (See Step 5 in the next section for more on this operation.) We have $(J_U)_{/F_U}\simeq \oplus_i \mathrm{Jac}(V_{S_i})$.

Via $\iota_1: B^{\times}_+\rightarrow \GL_2(\R)_+$, we have an action of $B^{\times}_+$ on the upper half plane $\HH$ by fractional linear transformations. Letting $\Gamma_i:=B^{\times}_+\cap S_i$, there is an isomorphism $V_{S_i}(\C)\simeq \Gamma_i\backslash\HH$ of compact Riemann surfaces. (This is true at least when $d>1$. For $d=1$ we would have to complete $\HH$ by the addition of cusps.) We also have $J_U(\C)\simeq \prod_{i=1}^h\mathrm{Jac}(V_{S_i})(\C)$.

Given a non-zero ideal $\mm$ of $\OOO_F$, we define 
$$U(\mm):=\{a\in B^{\times}_f: \ord_v(a_v-1)\geq\ord_v(\mm)\,\,\,\forall\text{ finite }v\}$$
$$=\{a\in B^{\times}_f: \iota_v(a_v)\equiv\begin{pmatrix} 1 & 0\\0 & 1\end{pmatrix}\pmod{\mm}\,\,\,\forall\text{ finite }v\}.$$
Any open compact subgroup of $B^{\times}_f$ contains $U(\mm)$ for some $\mm$.
Examples are 
$$U_{1,1}(\mm):=\{a\in B^{\times}_f: \iota_v(a_v)\equiv\begin{pmatrix} 1 & *\\0 & 1\end{pmatrix}\pmod{\mm}\,\,\,\forall\text{ finite }v\},$$
$$U_1(\mm):=\{a\in B^{\times}_f: \iota_v(a_v)\equiv\begin{pmatrix} * & *\\0 & 1\end{pmatrix}\pmod{\mm}\,\,\,\forall\text{ finite }v\},$$
and
$$U_0(\mm):=\{a\in B^{\times}_f: \iota_v(a_v)\equiv\begin{pmatrix} * & *\\0 & *\end{pmatrix}\pmod{\mm}\,\,\,\forall\text{ finite }v\}.$$
Note that $\nu(U_0(\mm))=\prod_{v\text{ finite}}\OOO_v^{\times}$, so $\A^{\times}/F^{\times}\nu(S(U_0(\mm)))$ is the narrow class group of $F$. (An equivalent description is as the group of fractional ideals modulo principal ideals generated by totally positive elements.) Thus if $F$ has narrow class number $1$ then $F_U=F$ and $X_U$ is geometrically connected. We let $X_0(\mm):=X_{U_0(\mm)}$, $X_1(\mm):=X_{U_1(\mm)}$, $X_{1,1}(\mm):=X_{U_{1,1}(\mm)}$, $J_0(\mm):=J_{U_0(\mm)}$, $J_1(\mm):=J_{U_1(\mm)}$ and $J_{1,1}(\mm):=J_{U_{1,1}(\mm)}$. 

Now let $E/F$ be an elliptic curve, of conductor $\n$. We say that $E$ is {\em geometrically modular} if there exists a surjective morphism $\pi : J_0(\n)\rightarrow E$, defined over $F$. We say that $E$ has a {\em hyperbolic uniformisation of arithmetic type} if for some open compact subgroup $U$ of $B^{\times}_f$, letting $\Gamma:=B^{\times}_+\cap S(U)$, there is a surjective holomorphic map
of Riemann surfaces $\theta: \Gamma\backslash\HH\rightarrow E(\C)$. Without loss of generality, we may assume that $U=U(\mm)$, for some non-zero ideal $\mm$ of $\OOO_F$. Then we have a surjective holomorphic map $\theta: \Gamma(\mm)\backslash\HH\rightarrow E(\C)$ of Riemann surfaces, where $\Gamma(\mm):=B^{\times}_+\cap S(U(\mm))$. (We define also $\Gamma_{1,1}(\mm)$, $\Gamma_1(\mm)$ and $\Gamma_0(\mm)$ in the obvious way.)

\begin{thm}\label{main} Let $F$ be a totally real field of odd degree, and $E/F$ an elliptic curve. If $E$ has a hyperbolic uniformisation of arithmetic type then $E$ is geometrically modular. 
\end{thm}

In the case that $F=\Q$, we have $B=\GL_2$, and the Shimura curves become modular curves. A proof in this case was sketched in the ``Technical Appendix'' to Mazur's expository article \cite{maz}, and expanded upon in more detail by Banwait \cite{bb}. Our proof of the more general case will be organised using steps largely analogous to those followed by Banwait. This follows in the next section. In Remark \ref{modular} we return to examine the relation between geometric modularity and other notions of modularity involving automorphic representations of $\GL_2(\A_F)$, Galois representations and numbers of points mod $\p$.

Aside from requiring the introduction of more background on Shimura curves, working in the generality that $F$ may not be $\Q$ causes a few difficulties of substance that must be dealt with. This is especially so in Step 4, including knock-on effects from the modification of Step 2, where we must use $U_{1,1}(\mm^2)$, not $U_1(\mm^2)$. Even if one was only interested in $F=\Q$, there are places where we take the opportunity to provide additional detail or clarification, for example in Steps 9 and 10.

\section{Proof of Theorem \ref{main}.}

\subsection*{Step 1. Reduce to the non-CM case.} 
In the case $F=\Q$, Shimura showed that an elliptic curve with complex multiplication is geometrically modular \cite{shimura1971elliptic}. More generally, for any totally real $F$, suppose that the elliptic curve $E/F$ has complex multiplication, so $\End_{\Qbar}(E)\otimes\Q\simeq K$, an imaginary quadratic field. By a theorem of Deuring \cite[Ch. II, Theorem 10.5(b)]{sil2}, there exists a Hecke Grossencharacter $\psi_{E/F'}:\A_{F'}^{\times}/{F'}^{\times}\rightarrow K^{\times}$ (where $F':=FK$) such that
$L(E/F,s)=L(s,\psi_{E/F'})$. This in turn is $L(s, \Ind_{F'}^F(\psi_{E/F'}))$, where $\Ind_{F'}^F(\psi_{E/F'})$ is a cuspidal automorphic representation of $\GL_2(\A_F)$, the automorphic induction of $\psi_{E/F'}$. Thus $E/F$ is modular in the sense of Remark \ref{modular} below, hence geometrically modular, as explained there. 

\subsection*{Step 2. Replace $U(\mm)$ with $U_{1,1}(\mm^2)$.}

Replacing $\mm$ by a smaller ideal if necessary, we may assume without loss of generality that $\mm=(M)$ is principal, with non-zero $M\in \OOO_F$ and $\tau_1(M)>0$. Then $M$ satisfies the hypotheses of \cite[Theorem 28.6.1]{jv} (Eichler's theorem on norms), so there exists $\alpha\in B_+^{\times}$ such that $\nu(\alpha)=M$. This element $\alpha$ may be diagonally embedded in $B^{\times}_f$. The element $\begin{pmatrix} M & 0\\0 & 1\end{pmatrix}$ of $\GL_2(F)$ is diagonally embedded in $\GL_2(\A_f)\simeq B^{\times}_f$. Letting $\beta:=\alpha^{-1}\begin{pmatrix} M & 0\\0 & 1\end{pmatrix}\in B^{\times}_f$, we have $\beta\in B^1_f:=\{b\in B^{\times}_f:\,\nu(b)=1\}$. Letting $B^1:=\{b\in B^{\times}:\,\nu(b)=1\}$, the Strong Approximation theorem \cite[Theorem 28.5.3]{jv} says that $B^1$ is dense in $B^1_f$. Hence there exists $\delta'\in B^1$ such that $\iota_v(\delta')\in\GL_2(\OOO_v)$ for all finite $v\nmid\mm$ and $\iota_v(\delta')$ is $v$-adically very close to $\beta_v$ for all finite $v\mid\mm$.

Now let $\delta=\alpha\delta'$ so that $\delta\in B_+^{\times}$, with $\iota_v(\delta)\in\GL_2(\OOO_v)$ for all finite $v\nmid\mm$ and $\iota_v(\delta)$ $v$-adically very close to $\begin{pmatrix} M & 0\\0 & 1\end{pmatrix}$ for all finite $v\mid\mm$. Let  $\Gamma=\Gamma(\mm)$ and $\tilde{\Gamma}:=\delta^{-1}\Gamma\delta$. If $\gamma'\in\Gamma_{1,1}(\mm^2)$ and $v\mid\mm$ then $\iota_v(\gamma')=\begin{pmatrix}a&b\\c&d\end{pmatrix}\equiv\begin{pmatrix} 1 & *\\0 & 1\end{pmatrix}\pmod{\mm^2}$, then
$\iota_v(\delta\gamma'\delta^{-1})$ is $v$-adically close to $\begin{pmatrix} a & bM\\c/M & d\end{pmatrix}$, hence in $\Gamma=\Gamma(\mm)$. We conclude that $\tilde{\Gamma}\supseteq\Gamma_{1,1}(\mm^2)$. 

If $\gamma\in\Gamma$ then for any $z\in\HH$ we have $\delta(\delta^{-1}\gamma\delta)z=\gamma(\delta z)$, hence $z\mapsto \delta z$ gives a well-defined holomorphic bijection from $\tilde{\Gamma}\backslash\HH$ to $\Gamma\backslash\HH$. We now have holomorphic surjections
$$\Gamma_{1,1}(\mm^2)\backslash\HH\rightarrow\tilde{\Gamma}\backslash\HH\rightarrow\Gamma\backslash\HH\rightarrow E(\C).$$
Re-labelling $\mm^2$ as $\mm$, since it is just a special case of a non-zero ideal of $\OOO_F$, we have 
$$\theta_1: \Gamma_{1,1}(\mm)\backslash\HH\rightarrow E(\C).$$

\subsection*{Step 3. Specialise from $\C$ to $\Qbar$ to a finite Galois extension $L/F$.} 

Recall the notation $S=S(U)=B^{\times}_{\infty +}U\leq B^{\times}_{\A +}$, where for us $U=U_{1,1}(\mm)$. The Shimura curve $X_U/F$ has geometric components $V_{S_i}/F_U$, where $S_i=x_iSx_i^{-1}$, with $B^{\times}_{\A +}= \cup_{i=1}^h B^{\times}_+ x_iS$ and $x_1=1$. Using the fact that a holomorphic map between complete nonsingular algebraic curves is necessarily an algebraic morphism (Riemann's existence theorem), we now have a surjective morphism of algebraic curves over $\C$, 
$\theta_1: V_S\rightarrow E$. Using \cite[\S 3.2.2]{bp}, as in \cite[\S2]{bb}, this can be replaced by a morphism defined over $\Qbar$. Roughly speaking, the curves and morphism are defined over a finitely generated subfield of $\C$. One may view this purely algebraically, then plug in algebraic values for algebraically independent transcendental generators.
In fact, clearly $\Qbar$ can be replaced by some finite Galois extension $L$ of $F$, which we had better choose to contain the abelian extension $F_U$ over which all the components $V_{S_i}$ are defined.

Recall that the Artin symbol $[\nu(x_i^{-1}), F]\in\Gal(F_U/F)$ takes $V_{S}$ to $V_{S_i}$. If we let $\tilde{\sigma_i}\in\Gal(L/F)$ be a lifting of $\sigma_i:=[\nu(x_i^{-1}), F]$ then we may conjugate by $\tilde{\sigma_i}$ the morphism $\theta_1: V(S)\rightarrow E$, defined over $L$, to get $\theta_i: V(S_i)\rightarrow E$, also defined over $L$. This is dependent on the choice of lifting, but fixing a choice we may view the collection of morphisms $\theta_i$ as a single morphism $\theta: X_{1,1}(\mm)\rightarrow E$, defined over $L$. This is because $X_{1,1}(\mm)_L$ is a disjoint union of the components $V_{S_i}$, each considered over $L$, and we just need some way of producing a morphism over $L$ from each such component to $E$.  

\subsection*{Step 4. Pass to the jacobian $J_{1,1}(\mm)$ and decompose.}

By the universal property of jacobians \cite[Proposition 6.4]{miljac}, applied separately to each geometrical component, from the morphism $\theta: X_{1,1}(\mm)\rightarrow E$, defined over $L$, we get a morphism of abelian varieties $\theta_*: J_{1,1}(\mm)\rightarrow E$, also defined over $L$. On the level of divisor classes, this is push-forward by $\theta$.

For each prime ideal $\p$ of $\OOO_F$ we have a Hecke correspondence $T_{\p}$ on $X_{1,1}(\mm)$, defined over $F$. (A correspondence from $A$ to $B$ is a certain subvariety of the product $A\times B$, generalising the graph of a morphism. In this case $A=B=X_{1,1}(\mm)$.) Over $F_U$ it is a disjoint union of correspondences between pairs of components $V_{S_i}$ conjugate by $\sigma(x):=[\mathfrak{w}_{\p}^{-1}, F]\in\Gal(F_U/F)$, where $x\in B^{\times}_{\A +}$ is such that $\iota_{\p}(x_{\p})=\left(\begin{pmatrix}\mathfrak{w}_{\p} & 0\\0 & 1\end{pmatrix}\right)$ (and $\mathfrak{w}_{\p}\in F$ is a uniformiser at $\p$), while $x_v=1$ at all other places. In the notation of \cite[Theorem 3.4.2]{oh} these are the correspondences $X_{S_iS_i}(x)$ from $V_{S_i}$ to $V_{S_i}^{\sigma(x)}$. See also \cite[beginning of \S 8]{kl}.

For each $\p\nmid\mm$ we also have a correspondence $\langle\p\rangle$, defined over $F$. Over $F_U$ it is a disjoint union of correspondences between pairs $V_{S_i}$ and $V_{S_i}^{\sigma(\mathfrak{w}_{\p})}$, $J_{S_iS_i}(x')$ in the notation of \cite[Theorem 2.3.1]{oh}, \cite[(2), 2.5]{shi70}, where $x'$ is obtained from $x$ upon replacing $x_{\p}$ by $\mathfrak{w}_{\p}$. This $\langle\p\rangle$ depends only on the class of $\p$ (or rather of the element $y\in\A^{\times}$ that is $\mathfrak{w}_{\p}$ at $\p$,  $1$ at other places) in the ray-class group $C_{\mm +}:=\A^{\times}/F^{\times}F^{\times}_{\infty +}\A_f^{1,\mm}$, where $\A_f^{1,\mm}:=\{w\in \A_f^{\times}:\, \ord_v(w_v-1)\geq\ord_v(\mm)\,\,\forall v\}$. 

Consider a class $c\in U_0(\mm)/U_{1,1}(\mm)\simeq ((\OOO_F/\mm)^{\times})^2$, represented by an element $g_c\in B^{\times}_f$. Since $g_cU_{1,1}(\mm)g_c^{-1}\subseteq U_{1,1}(\mm)$, we have a morphism defined over $F$, $g_c: X_{1,1}(\mm)\rightarrow X_{1,1}(\mm)$. (Such morphisms associated with elements of $B^{\times}_f$ are recalled just after (1.1.2) in \cite{zh}.) Then we have an endomorphism defined over $F$, $g_c^*: J_{1,1}(\mm)\rightarrow J_{1,1}(\mm)$, depending only on the class $c$. 

Let $\TT'$ be the subring of $\End_F(J_{1,1}(\mm))$ generated (over $\ZZ$) by the $T_{\p}$ and $\langle \p\rangle$ for $\pp\nmid\mm$, and by the $g_c^*$ for all $c\in U_0(\mm)/U_{1,1}(\mm)$. Given a height-zero prime ideal $I$ in $\TT'$, we define $J_I$ (an abelian variety defined over $F$) to be the connected component of the identity in the kernel of $I$ on $J_{1,1}(\mm)$. To get a feel for such $I$ and $J_I$, we need to consider cotangent spaces and automorphic forms.

There are natural pullback actions, defined over $F$, of each $T_{\p}$, $\langle\p\rangle$ and $g_c$ on the $F$-vector space of regular differentials $H^0(X_{1,1}(\mm), \Omega^1)$, which is the cotangent space to $J_{1,1}(\mm)$ at the identity. Over $\C$, 
$H^0(X_{1,1}(\mm)_{\C}, \Omega^1)\simeq \oplus_{i=1}^h H^0(\Gamma_i\backslash\HH,\Omega^1)$, where $B^{\times}_{\A +}$ is a disjoint union $\cup_{i=1}^h B^{\times}_+ x_iS(U_{1,1}(\mm))$. An element $\omega=(f_i(\tau)\,d\tau)_{i=1}^h$ of this space may be identified with a function $\phi_{\omega}$ on $B_{\A }^{\times}$ in the following way. On an element $g=\gamma x_ig_1g_2\ldots g_d u$, with $1\leq i\leq h$, $\gamma\in B^{\times}$, $g_j\in B_{\infty_j +}$ $\left(\iota_{\infty_1}(g_1)=\begin{pmatrix}a & b\\c & d\end{pmatrix}\right)$ and $u\in U_{1,1}(\mm)$, it takes the value $f_i(\tau)(c\tau+d)^2$, with $\tau=g_1(\sqrt{-1})\in\HH$. This function $\phi_{\omega}$ is well-defined, and satisfies $\phi_{\omega}(\gamma gu)=\phi_{\omega}(g)\,\,\,\forall g\in B^{\times}_{\A}, \gamma\in B^{\times}$ and $u\in U_{1,1}(\mm)$. Under the right action of $U_0(\mm)/U_{1,1}(\mm)\simeq ((\OOO_F/\mm)^{\times})^2$, suppose that it transforms by a character $(\chi_1,\chi_2)$, i.e. $\begin{pmatrix}a & b\\c & d\end{pmatrix}\mapsto \chi_1(a)\chi_2(d)$. Denoting also by $\chi_1$ and $\chi_2$ the lifts to $\A^{\times}$, we have $\phi_{\omega}\in L^2(B^{\times}\backslash B^{\times}_{\A},\chi_1\chi_2)$, the ``$\chi_1\chi_2$'' denoting that $\phi_{\omega}(zg)=\chi_1\chi_2(z)\phi(g)\,\,\,\forall z\in\A^{\times}$. 

In classical terms, $H^0(X_{1,1}(\mm)_{\CC}, \Omega^1)$ is a direct sum over $1\leq i\leq h$ of spaces of quaternionic modular forms of weight $2$, holomorphic functions on $\HH$ transforming appropriately under the $\Gamma_i$. The actions of $T_{\p}$ and $\langle\p\rangle$ become those of 
$\mathfrak{T}(SxS)$ and $\mathfrak{T}(Sx'S)$, in the notation of \cite[\S 1.5]{oh}. (Note that in general these operators permute the summands.)
The space $L^2(B^{\times}\backslash B^{\times}_{\A},\chi_1\chi_2)$, with the right-translation action of $B^{\times}_{\A}$, decomposes as a direct sum of irreducible automorphic representations $\pi\simeq\otimes_v \pi_v$. Because it comes from a quaternionic modular form of parallel weight $2$, the vector $\phi_{\omega}=:w=\otimes_v w_v$ lies in the subsum of those $\pi$ such that $\pi_{\infty_1}$ is the discrete series representation $D_2$ of $\GL_2(\R)$ (with $w_{\infty_1}$ in the $1$-dimensional space of holomorphic vectors in the lowest $\mathrm{SO}_2$-type) and $\pi_{\infty_j}$ is trivial for $2\leq j\leq d$. For finite $v\nmid\mm$, $w_v$ is in the $1$-dimensional space of $\GL_2(\OOO_v)$-fixed vectors. 

The $\langle\p\rangle$ and $T_{\p}$ (for $\p\nmid\mm$) and the $[g_c]\in U_0(\mm)/U_{1,1}(\mm)$ are simultaneously diagonalisable on $H^0(X_{1,1}(\mm)_\C, \Omega^1)$, with $U_0(\mm)/U_{1,1}(\mm)\simeq ((\OOO_F/\mm)^{\times})^2$ acting on an eigenvector via a character $(\chi_1,\chi_2)$, and $\langle\p\rangle$ as $\chi_1\chi_2(\p)$. For each eigenvector $\mathrm{v}$, there is a homomorphism $\theta_{\mathrm{v}}:\TT'\rightarrow \CC$, taking each operator to the eigenvalue by which it acts on $\mathrm{v}$. The kernel of $\theta_{\mathrm{v}}$ is one of the height-zero prime ideals $I$, all of which arise this way. Eigenvectors whose systems of eigenvalues are Galois conjugate give rise to the same  $I$. So the cotangent space of $J_I$ is a direct sum of Galois conjugacy classes of eigenspaces. The jacobian $J_{1,1}(\mm)$ is isogenous to a direct sum $\oplus_I J_I$. We must now reorganise this decomposition.  

First we look at the subsum of those $J_I$ for which $\chi_1=1$. This is isogenous to $J_1(\mm)$, with cotangent space $H^0(X_1(\mm), \Omega^1)$. Since $B^{\times}_f\simeq\GL_2(\A_f)$ there is a theory of newforms just as for $\GL_2$ (cf. \cite[\S 3.1]{zh}), and on the new subspace of $H^0(X_1(\mm)_{\C}, \Omega^1)$ the $T_{\p}$ for $\p\mid\mm$ are also diagonalisable. The simultaneous new eigenspaces are $1$-dimensional. This is as in \cite[Theorem 3.2.1(2)]{zh}, but easier, because $B^{\times}_f\simeq\GL_2(\A_f)$ puts us always in his Case 1 (though our $\chi=\chi_2$ is not necessarily trivial). It depends on the following two facts, in addition to the observations about $1$-dimensionality in the previous-but-one paragraph.
\begin{enumerate}
\item The multiplicity of $\pi$ in $L^2(B^{\times}\backslash B^{\times}_{\A},\chi_2)$ is $1$. This is \cite[Theorem 10.10]{G}, and is a consequence of the trace formula proof of the Jacquet-Langlands correspondence, and the analogous theorem for $\GL_2$.
\item If $\mathfrak{M}$ is maximal such that $\pi_f$ has $U_1(\mathfrak{M})$-fixed vectors then the dimension of the fixed subspace is $1$.
\end{enumerate}

It follows, as in \cite[Lemma 3.4.5(1)]{zh}, that $J_1(\mm)$ is isogenous over $F$ to a direct sum $\bigoplus_{[\Phi]}J_{\Phi}^{m_{\Phi}}$, where $[\Phi]$ runs through the Galois equivalence classes of newforms at levels $\mathfrak{M}_{\Phi}\mid\mm$ and $m_{\Phi}$ is the number of divisors of $\frac{\mm}{\mathfrak{M}_{\Phi}}$. Each $J_{\Phi}$ is a simple abelian variety over $F$, with endomorphism ring containing a copy of $\OOO_{\Phi}$, the subring of $\C$ generated by the eigenvalues $a_{\p}(\Phi)$ and $\chi_2(\p)$ on $\Phi$ of $T_{\p}$ and $\langle\p\rangle$ (respectively), for $\p\nmid \mathfrak{M}_{\Phi}$, via the action of the Hecke correspondences on $J_1(\mathfrak{M}_{\Phi})$. It can be constructed as the largest abelian subvariety of $J_1(\mathfrak{M}_{\Phi})$ killed by the kernel of the homomorphism $\TT\rightarrow\C$ sending each operator to its eigenvalue, where $\TT$ is the subring of $\End_F(J_1(\mathfrak{M}_{\Phi}))$ generated by the $T_{\p}$ and $\langle\p\rangle$, for $\p\nmid \mathfrak{M}_{\Phi}$.  It is simple because its endomorphism algebra contains the fraction field of $\OOO_{\Phi}$, whose degree over $\Q$ is the same as the dimension of $J_{\Phi}$, which is the dimension of its cotangent space, the size of a Galois equivalence class of newforms.

Zhang's $J_{\Phi}$ is like our $J_{\Phi}^{m_{\Phi}}$, because he does this construction inside $J_1(\mm)$ rather than $J_1(\mathfrak{M_{\Phi}})$. In fact it is one of our $J_I$. Let's call it $J_{I(\Phi)}$. The isogeny from $J_{\Phi}^{m_{\Phi}}$ to $J_{I(\Phi)}$ is a sum of restrictions to $J_{\Phi}$ of $g_{\mathfrak{d}}^*: J_1(\mathfrak{M}_{\Phi})\rightarrow J_1(\mm)$, where $\mathfrak{d}$ runs over divisors of $\frac{\mm}{\mathfrak{M}_{\Phi}}$, $g_{\mathfrak{d}}\in B^{\times}_f$ has component $\iota_v^{-1}\left(\begin{pmatrix} \mathfrak{w}_v^{\ord_v(\mathfrak{d})} & 0\\0 & 1\end{pmatrix}\right)$ at each finite place $v$, and we have a morphism defined over $F$, $g_{\mathfrak{d}}: X_1(\mm)\rightarrow X_1(\mathfrak{M}_{\Phi})$, since $g_{\mathfrak{d}}U_1(\mm)g_{\mathfrak{d}}^{-1}\subseteq U_1(\mathfrak{M}_{\Phi})$.  For the case $F=\QQ$ of this decomposition of the modular Jacobian, see \cite[Proposition 2.3]{rib}, or \cite[Theorem 6.6.6]{ds} for $J_0(N)$.

We claim that the subsums of $J_I$ for each non-trivial $\chi_1$ may be described in exactly the same way according to Galois conjugacy classes of newforms $\Phi$, but with values of $\chi_1$ also thrown in to $\OOO_{\Phi}$. This is because of a certain operator 
$$\theta_{\chi_1}: H^0(X_{1,1}(\mm)_{\C}, \Omega^1)^{(\chi_1,\chi_2)}\rightarrow H^0(X_{1,1}(\mm)_{\C}, \Omega^1)^{(1,\chi_1\chi_2)},$$
where the superscripts indicate that we are considering subspaces on which $U_0(\mm)/U_{1,1}(\mm)\simeq ((\OOO_F/\mm)^{\times})^2$ acts by these characters. Viewing these spaces of automorphic forms as functions $f$ on $B^{\times}_{\A}$, and letting $f^{\chi_1}:=\theta_{\chi_1}(f)$, it is defined by
\begin{equation}\label{twist} f^{\chi_1}(g):=\sum_{a\,(mod\,\mm)}\chi_1(a) f\left(g\begin{pmatrix} \mm & a\\0 & \mm\end{pmatrix}\right),\end{equation} where $\begin{pmatrix} \mm & a\\0 & \mm\end{pmatrix}$ is short for the element of $B^{\times}_f$ that is $\iota_v^{-1}\begin{pmatrix}\mathfrak{w}_v^{\ord_v(\mm)} & a\\0 & \mathfrak{w}_v^{\ord_v(\mm)}\end{pmatrix}$ at finite $v\mid\mm$, $1$ elsewhere.

One checks that for $t=\begin{pmatrix}t_1 & 0\\0 & t_2\end{pmatrix}\in B^{\times}_f\simeq\GL_2(\A_f)$,
$$f^{\chi_1}(gt)=\sum_{a\,(mod\,\mm)}\chi_1(a) f\left(gt\begin{pmatrix} \mm & a\\0 & \mm\end{pmatrix}\right)$$
$$=\chi_1(t_2t_1^{-1})\sum_{a\,(mod\,\mm)}\chi_1(t_1at_2^{-1})f\left(g\begin{pmatrix}\mm & t_1at_2^{-1}\\0 & \mm\end{pmatrix}t\right)=(\chi_1,\chi_2)(t)\chi_1(t_2t_1^{-1})f^{\chi_1}(g).$$
Hence $f^{\chi_1}(gt)=(1,\chi_1\chi_2)(t) f^{\chi_1}(g)$, as required. 

There is a similar operator going in the opposite direction, providing an inverse to $\theta_{\chi_1}$, which is therefore an isomorphism. Thus, the decomposition of $H^0(X_{1,1}(\mm)_{\C}, \Omega^1)^{(1,\chi_1\chi_2)}\simeq H^0(X_1(\mm)_{\C}, \Omega^1)^{\chi_1\chi_2}$ according to newforms is mirrored in $H^0(X_{1,1}(\mm)_{\C}, \Omega^1)^{(\chi_1,\chi_2)}$. In fact we are looking at different vectors in the same automorphic representations, with the same central character $\chi_1\chi_2$. The isogeny decomposition $$J_{1,1}(\mm)\sim_F\bigoplus_{[\Phi]}(J_{\Phi/L})^{m_{\Phi}}$$ follows, where now the newforms $\Phi$ range over both trivial and non-trivial $\chi_1$.

So we now have a morphism of abelian varieties $$\theta_*: \bigoplus_{[\Phi]}(J_{\Phi/L})^{m_{\Phi}}\rightarrow E_{/L}.$$

\subsection*{Step 5. Take the Weil restriction of scalars down to $F$.} 

The Weil restriction of scalars produces from an abelian variety $A$ over $L$ an abelian variety $\Res^L_F(A)$ over $F$, of $r:=[L:F]$ times the dimension, with $\Res^L_F(A)_{/L}\simeq \oplus_{\sigma\in\Gal(L/F)}A^{\sigma}$, cf. \cite[p.178, l.1]{milne1972arithmetic}. In particular, if $A$ is already defined over $F$ then $\Res^L_F(A)_{/L}\simeq A^r$. Given a homomorphism defined over $L$, $f: A\rightarrow B$, between abelian varieties defined over $F$, we have a map from $\Res^L_F(A)_{/L}\rightarrow B$ such that $(a_1,\ldots,a_r)\mapsto f(a_1)+\ldots +f(a_r)$. We may view this as a $\Res^L_F(A)_{/L}$-valued point of $B_{/L}$. Using a key property of $\Res^L_F$, i.e. \cite[(4.1)]{mazur2007twisting}, \cite[\S 1, l.5]{milne1972arithmetic} this is equivalent to a $\Res^L_F(A)$-valued point of $\Res^L_F(B)$. This is a morphism over $F$ from $\Res^L_F(A)$ to $\Res^L_F(B)$, necessarily a homomorphism of abelian varieties, by \cite[Corollary 2.2]{milab}. 

Given an abelian variety $C$ over $F$, and a $\Z[G]$-module $M$, free of some rank $R$ over $\Z$, where $G:=\Gal(L/F)$, there is also an abelian variety $C\otimes M$ over $F$, a twisted form of $C^R$, so $(C\otimes M)_{/L}\simeq (C_{/L})^R$, as in \cite[Definition 1.1]{mazur2007twisting} or \cite[\S 2]{milne1972arithmetic}. (Note that for any $\sigma\in\Gal(L/F)$, $C^{\sigma}\simeq C$, since $C$ is defined over $F$.) Moreover, by \cite[Proposition 4.1]{mazur2007twisting} or \cite[above Proposition 6]{milne1972arithmetic}, we have an isomorphism $C\otimes \Z[G]\simeq \Res^L_F(C_{/L})$. Applying $\Res^L_F$ to $\theta_*$ above, we get then
\begin{equation}\label{eq:main}
\theta_*: \bigoplus J_{\Phi}^{m_{\Phi}}\otimes\Z[G]\rightarrow E\otimes\Z[G].\end{equation}

\subsection*{Step 6. Identify, up to $F$-isogeny, some $J_{\Phi}$ with a twist of $E$ by an irreducible representation $M$ of $G$.} 
The $\Q[G]$-module $\Q[G]$ decomposes as a direct sum of irreducibles: $\Q[G]=\oplus V_i$. Note that these representations of $G$ are irreducible over $\Q$, not necessarily over $\C$. Now for each $i$ define $M_i:=V_i\cap\Z[G]$. Then $\oplus M_i$ is a $\Z[G]$-module of finite index in $\Z[G]$. It follows that given an abelian variety $C$ over $F$, $C\otimes\Z[G]$ is isogenous to $\oplus_i (C\otimes M_i)$.

On the left-hand side of (\ref{eq:main}), each $J_{\Phi}$ arises as a simple factor, because $\Q[G]$ contains the trivial representation. To prove that at least one $J_{\Phi}$ must be identified with an isogeny factor of the right hand side, the only thing that could go wrong is if the map all the way from $J_{1,1}(\mm)$ via $J_{1,1}(\mm)\otimes \Z[G]$ to $E \otimes \Z[G]$ were zero. To get a contradiction, base change to $L$, and we have $J_{1,1}(\mm)_{/L}$ to $(J_{1,1}(\mm)_{/L})^r$ to $(E_{/L})^r$, where $r=[L:F]$. The first map is the diagonal embedding, so by surjectivity of $\theta_*: J_{1,1}(\mm)_{/L} \rightarrow E_{/L}$ the composite cannot be $0$. 

Therefore some $J_{\Phi}$ arises as an isogeny factor of the right hand side. On the other hand, the right hand side is isogenous to $\oplus_i (E\otimes M_i)$, as above. \emph{A priori} these twists need not be simple over $F$; however here they are, since $E$ does not have CM \cite[Proposition 2.6]{bartel2013simplicity}. We therefore obtain, for some $J_{\Phi}$, an irreducible $\Z[G]$-module $M_{\Phi}$ such that
\[ J_{\Phi} \sim_F E \otimes_\Q M_{\Phi}. \]
At this point, one wants to ``untwist $E$'' to obtain
\[ J_{\Phi} \otimes_\Q M_{\Phi}^{-1} \sim_F E, \]
and then to relate the left-hand side here with the abelian variety associated to a newform. Doing this rigorously is essentially the content of the rest of the proof. 

\medskip In the sequel, we fix one such $\Phi$ as above, and for ease of notation, we replace $M_{\Phi}$ by $M$.

\subsection*{Step 7. Show that $G$ acts on $M$ via a character $\chi'$.} Since the dimension of $J_{\Phi}$ is $[K_{\Phi}:\Q]$, as in \cite[Lemma 3.4.5(3)]{zh}, and the dimension of $E \otimes M$ is $\rank_{\Z}(M)$, we obtain 
\[ [K_{\Phi}:\Q] = \rank_{\Z}(M)=\dim_{\Q}(V),\]
where $V:=M\otimes_{\Z}\Q$. Here $K_{\Phi}$ is the fraction field of $\OOO_{\Phi}$. Via Hecke correspondences we have an embedding 
\[ K_{\Phi}\hookrightarrow (\End_F(J_{\Phi}))\otimes_{\Z}\Q = (\End_L(E)\otimes_{\Z}\End_{\Q}(V))^G=(\End_{\Q}(V))^G, \]
with the penultimate equality coming from \cite[Lemma 2.1]{bartel2013simplicity}. This establishes $V$ as a $K_{\Phi}$-vector space of dimension $1$, the $K_{\Phi}$-action commuting with that of $G$. Since $\dim_{K_{\Phi}}(V) = 1$, we obtain that the $G$-action may be expressed via a character
\[ \chi': G \to K_{\Phi}^\times. \]

\subsection*{Step 8. Identify $\chi'$ as a Dirichlet character.} Since $K_{\Phi}^\times$ is abelian, $\chi'$ factors through the abelianization $G^{\mathrm{ab}}=\Gal(L^{\mathrm{ab}}/F)$ of $G$, where $L^{\mathrm{ab}}/F$ is the maximal abelian subextension of $L/F$.  For some suitable ideal $\mm'$ of $\OOO_F$, $L^{\mathrm{ab}}$ is a subfield of the ray class field $H_{\mm' +}$ (conductor $\mm'\prod_{i=1}^{d}\infty_i$). Thus by class field theory $\chi'$ may be viewed as the restriction of a character (also denoted $\chi'$) of the ray class group $C_{\mm' +}$. Here $H_{m' +}$ denotes the narrow ray class field of conductor
$m'\prod_{i=1}^d \infty_i$, and the identification of
$\Gal(H_{m'}^+/F)^\vee$ with the character group of the narrow ray class group
$C_{m' +}$ is the standard reciprocity isomorphism from global class field
theory; see for example \cite[Ch. VI]{ne}.

\subsection*{Step 9. Show $a_{\p}(\Phi) = a_{\p}(E)$ for enough $\p$.} Here $a_{\p}(\Phi)$ is the eigenvalue of $T_{\p}$ on the newform $\Phi$, and (for $\p\nmid\n$) $a_{\p}(E)$ is such that $\#\tilde{E}(\FF_{\p})=1+\Nm(\p)-a_{\p}(E)$, where $\tilde{E}/\FF_{\p}$ is the reduction. Our treatment is analogous to what Diamond and Shurman \cite[\S 8.8]{ds} do in the case $F=\Q$, but whereas they are dealing with $X_0(M)$, for us it is $X_{1,1}(\mm)$, so the operators $\langle\p\rangle$ are not trivial. 

We need a Shimura curve version of the Eichler-Shimura congruence relation. This was proved by Shimura \cite[Theorem 2.23]{shi70} for all but finitely many $\p$, without specifying precisely which $\p$. It was sharpened by Ohta \cite[Corollary 3.4.4]{oh}, who proved that it works for all $\p\nmid\mm$. (The additional condition that $B$ is split at $\p$ is automatic for us, since our $B$ is split at all finite places.) This congruence relation is an equality, of correspondences between reductions $\widetilde{V_{S_i}}$ and $\widetilde{V_{S_i}^{\sigma(x)}}$ , modulo $\frak{P}\mid\p$ in $F_U$. But putting it all together it becomes an equality
$$\widetilde{T_{\p}}=\phi_{\p}^t+\widetilde{\langle\p\rangle}\circ\phi_{\p}$$ of correspondences on $\widetilde{X_{1,1}(\mm)}$, where $\phi_{\p}$ is the $\Nm(\p)$-power Frobenius correspondence. 

Now suppose that $\p\nmid\mm\n\mm'$, so that $J_{\Phi}$ and $E$ have good reduction and $\chi'(\p)$ is defined. Suppose also that the class of $\p$ in $C_{\mm\mm' +}$ is trivial, so $\chi(\p)=1$ and $\langle\p\rangle_\ast$ acts trivially on $J_{\Phi}$. (This uses only triviality of the class of $\p$ in $C_{\mm +}$, but we shall see below where $\mm'$ is used.) This allows us to construct the following commutative diagram, allowing the lower rectangle to commute without us having to attempt to make sense of $\widetilde{\langle \p \rangle}_{\ast}$ on $\widetilde{E\otimes M}$. (Any such attempt is not spelled out in \cite{bb}.) Note that the commutativity of the top left rectangle follows from the above congruence relation.
\[
\begin{tikzcd}
J_{\Phi} \arrow{r}{a_{\p}(\Phi) - a_{\p}(E)} \arrow{d} &[6em] J_{\Phi} \arrow{r}{\alpha} \arrow{d} &[6em] E \otimes M \arrow{d}\\
\widetilde{J_{\Phi}} \arrow{r}{\phi_{\p}^\ast + \widetilde{\langle \p \rangle}_{\ast}\phi_{\p,\ast} - a_{\p}(E)} \arrow{d}{1} &[6em] \widetilde{J_{\Phi}} \arrow{r}{\tilde{\alpha}} &[6em] \widetilde{E \otimes M} \arrow{d}{1}\\
\widetilde{J_{\Phi}} \arrow{r}{\tilde{\alpha}} &[6em] \widetilde{E \otimes M} \arrow{r}{\phi_{\p}^\ast + \phi_{\p,\ast} - a_{\p}(E)}  &[6em] \widetilde{E \otimes M}.
\end{tikzcd}
\]
Here $\alpha$ is the isogeny from Step 6, and the top vertical maps are reduction mod $\p$.

If $a_{\p}(\Phi) \neq a_{\p}(E)$ then (using the fact that $J_{\Phi}$ is simple) the (composite of maps along the) top row is surjective, as is the vertical map $E \otimes M \to \widetilde{E \otimes M}$. But the bottom right horizontal map is zero, using the fact that on $\widetilde{E}$ multiplication by $a_{\p}(E)$ is the sum of the Frobenius morphism and its dual. We are also using the fact that over $\FF_{\p}$, $\widetilde{E \otimes M}\simeq \widetilde{E}^{[K_{\Phi}:\Q]}$, since $\p$ splits completely in $L^{\mathrm{ab}}\subseteq H_{\mm' +}$. Hence the middle row is zero, giving a contradiction because all the rectangles commute. We thus obtain, for all primes $\p$ as above, that
\[ a_{\p}(\Phi) = a_{\p}(E). \] 

\subsection*{Step 10. Conclude via Galois representations and the Isogeny theorem.}

For any prime number $\ell$ we have a $2$-dimensional $\ell$-adic representation $\rho_{\ell,E}$ of $\Gal(\Qbar/F)$ on $V_{\ell}(E):=\Q_{\ell}\otimes_{\Z_{\ell}}\varprojlim E[\ell^n]$. Similarly we have $V_{\ell}(J_{\Phi})$, which is a free module of rank $2$ over $K_{\Phi}\otimes_{\Q}\Q_{\ell}\simeq\oplus_{\lambda\mid\ell}K_{\Phi,\lambda}$, giving us a direct sum $V_{\ell}(J_{\Phi})\simeq\oplus_{\lambda\mid\ell}V_{\lambda}(J_{\Phi})$. We fix a choice of any $\lambda$, giving us a $2$-dimensional $K_{\Phi, \lambda}$-representation $\rho_{\lambda, J_{\Phi}}$ of $\Gal(\Qbar/F)$.

Then, for all primes  $\p\nmid\mm\n\mm'$ with the class of $\p$ in $C_{\mm\mm' +}$ trivial, and $\Frob_{\p}\in\Gal(\Qbar/F)$ a Frobenius element,
$$\Tr(\Frob_{\p}|V_{\ell}(E))=a_{\p}(E)=a_{\p}(\Phi)=\Tr(\Frob_{\p}|V_{\lambda}(J_{\Phi})).$$
(The last equality is a consequence of the congruence relation after reduction mod $\p$.)
To say that the class of $\p$ in $C_{\mm\mm' +}$ is trivial is equivalent to saying that $\pp$ splits completely in the ray class field $H_{\mm\mm' +}$. The divisors in $H_{\mm\mm' +}$ of such primes include all unramified primes in $H_{\mm\mm' +}$ of absolute degree $1$ (i.e. norm $p$). Hence the set of $\Frob_{\p}$ as above is dense in $\Gal(\Qbar/H_{\mm\mm' +})$, and it follows that $\rho_{\ell, E}$ and $\rho_{\lambda, J_{\Phi}}$ are isomorphic (over $K_{\Phi, \lambda}$) when restricted to $\Gal(\Qbar/H_{\mm\mm' +})$. (Note that this restriction is irreducible, by the main theorem of \cite{serreMG}, since $E$ does not have complex multiplication.)

By Frobenius reciprocity, if  $\chi''$ is a character of $\Gal(H_{\mm\mm' +}/F)$, and $L$ is a finite extension of $\Q_{\ell}$ sufficiently large to accommodate both $K_{\Phi,\lambda}$ and the values of all such $\chi''$, then 
$$\Hom_{L[\Gal(\Qbar/F)]}\left(\Ind_{\Gal(\Qbar/H_{\mm\mm' +})}^{\Gal(\Qbar/F)}(V_{\ell}(E)_L), V_{\ell}(E)_L\otimes\chi''\right)$$ $$\simeq \Hom_{L[\Gal(\Qbar/H_{\mm\mm' +})]}\left(V_{\ell}(E)_L, V_{\ell}(E)_L\otimes\chi''\right),$$
which is $1$-dimensional since these are isomorphic, irreducible restrictions. For dimension reasons the $\rho_{\ell}(E)_L\otimes\chi''$, as $\chi''$ runs through characters of $\Gal(H_{\mm\mm' +}/F)$, must be precisely the irreducible constituents of the induced representation. (Note that the assumption that $E$ has no CM implies that these twists are all non-isomorphic. If not, then the set of primes for which $a_{\pp}(E)=0$ would have positive density, contrary to the exercises in \cite[IV,2.2]{serreMG}.) But a similar application of Frobenius reciprocity shows that $\rho_{\lambda, J_{\Phi}}$ (over $L$) is also isomorphic to an irreducible constituent of $\Ind_{\Gal(\Qbar/H_{\mm\mm' +})}^{\Gal(\Qbar/F)}(V_{\ell}(E)_L)$.

Hence the representations $\rho_{\ell, E}\otimes\chi''$ and $\rho_{\lambda, J_{\Phi}}$ of $\Gal(\Qbar/F)$ are isomorphic, for some character $\chi''$ of $\Gal(H_{\mm\mm' +}/F)\simeq C_{\mm\mm' +}$. Replacing $\Phi$ by a newform in $\pi_{\Phi}\otimes(\chi''^{-1}\circ\nu)$, we may assume without loss of generality that $\chi''$ is trivial. Using the twisting operator $\theta_{\chi_1}$ from (\ref{twist}) (which does not change eigenvalues of $T_{\p}$ or $\langle \p\rangle$), we may assume also that $\chi_1=1$. So we have a newform $\Phi$ of level $U_1(\mm'')$, with $\mm''\mid(\mm\mm')^2$, and character $\chi_2=\psi$, say, with $\rho_{\ell, E}\simeq \rho_{\lambda, J_{\Phi}}$ (over $K_{\Phi, \lambda}$), as representations of $\Gal(\Qbar/F)$.

For any $\pp\nmid\mm''$ we have
\begin{equation}\label{traces} a_{\p}(\Phi)=\Tr(\rho_{\lambda, J_{\Phi}}(\Frob_{\pp}))=\Tr(\rho_{\ell, E}(\Frob_{\pp}))=a_{\pp}(E)\end{equation}
and $$\psi(\p)=\det(\rho_{\lambda, J_{\Phi}}(\Frob_{\p}))/\Nm(\p)=\det(\rho_{\ell, E}(\Frob_{\p}))/\Nm(\p)=1.$$
Hence $K_{\Phi}=\Q$ and $\psi$ is trivial, i.e. $J_{\Phi}$ is an elliptic curve and $\Phi$ is a newform for $U_0(\mm'')$.
By (\ref{traces}), $\rho_{\ell, E}$ and $\rho_{\lambda, J_{\Phi}}$ are isomorphic representations of $\Gal(\Qbar/F)$. Then by Faltings' Isogeny Theorem \cite{fa}, $E$ is isogenous to $J_{\Phi}$ over $F$, hence there is a surjective morphism of abelian varieties over $F$, $J_0(\mm'')\rightarrow E$. Comparing conductors of $\rho_{\ell, E}$ and $\rho_{\lambda, J_{\Phi}}$, $\n=\mm''$. (That the conductor of $\rho_{\lambda, J_{\Phi}}$ is $\mm''$ follows from a theorem of Carayol \cite{car}.) Hence $E/F$ is geometrically modular.

\begin{remar}\label{modular}
We could say that an elliptic curve $E$ defined over a totally real field $F$ is {\em modular} if there exists an irreducible, cuspidal, automorphic representation $\pi$ of $\GL_2(\A_F)$ such that $a_{\p}(E)=a_{\p}(\pi)$ for all but finitely many prime ideals $\p$ of $\OOO_F$. Here $\tilde{E}(\FF_{\p})=1+\Nm(\p)-a_{\p}(E)$ for a prime of good reduction, and $a_{\p}(\pi)$ is the eigenvalue of $T_{\p}$ on a $\GL_2(\OOO_{\p})$-fixed vector for $\pi_{\p}$, at a prime where $\pi_{\p}$ is spherical. (In classical terms, the $a_{\p}(\pi)$ are the Hecke eigenvalues for a cuspidal Hilbert modular form, an $h$-tuple of certain functions on $\HH^d$, where $h$ is the narrow class number of $F$ and $d=[F:\Q]$.) It would follow that the $\ell$-adic representations of $\Gal(\Qbar/F)$ attached to $E$ and (\cite{car}) to $\pi$ are isomorphic, since $a_{\p}(E)$ and $a_{\p}(\pi)$ are traces of $\Frob_{\p}$. Equivalently the $L$-functions $L(E,s)$ and $L(\pi, s)$ are the same. The primes of good reduction for $E$ will be the same as the spherical primes for $\pi$, and $\pi$ will have $U_0(\n)$-fixed vectors, unique up to scaling, where $\n$ is the conductor of $E$.

Under our assumption that $d$ is odd, and choosing the quaternion algebra $B/F$ to be split at all finite places and one real place, the Jacquet-Langlands correspondence \cite{jl} ensures that associated to $\pi$ is an automorphic representation of $B_{\A}^{\times}$, containing a $U_0(\n)$-newform $\Phi$, with the same Hecke eigenvalues as $\pi$, for $\p\nmid\n$. Then (\ref{traces}) holds, so the Galois representation $\rho_{\ell, \Phi}$ is the same as $\rho_{\ell, E}$, and by Faltings' Isogeny Theorem $E/F$ is geometrically modular. Conversely, if $E/F$ is geometrically modular then the Jacquet-Langlands correspondence provides a $\pi$ associated to $\Phi$, with respect to which $E/F$ is modular. 

Modularity of elliptic curves over totally real fields is known for $d=2$ \cite{fhs}, $d=3$ \cite{dns}, and for $d=4$ when $\sqrt{5}\notin F$ \cite{box}.
\end{remar}

\begin{remar}\label{deven} In attempting to remove the conditon that $d=[F:\Q]$ is odd, one would have to assume the existence of a finite place $\q$ of bad reduction for $E/F$, and let $B$ be the quaternion algebra over $F$ split at one real place and all finite places except $\q$. Locally at $\q$, the open compact subgroup $U_0(\q^r)$ might be replaced by the subgroup of $B^{\times}_{\q}$ generated by $U(\q^{r-1})$ and the group of units $\OOO_K^{\times}$ of the ring of integers of an embedded unramified quadratic extension $K$ of $F_{\q}$. The characters $(\chi_1,\chi_2)$ of $U_0(\mm)/U_{1,1}(\mm)$ would be replaced (locally at $\q$) by characters $\chi$ of $\OOO_K^{\times}$. One would then hope to exploit the theory of newforms developed in \cite[Propositions 3.3.1, 3.3.2]{zh} and by Gross \cite{Gr}. But the fact that they deal only with unramified characters $\chi$ of $K^{\times}$, restricting to the trivial character of $\OOO_K^{\times}$ (cf. \cite[(6.1)]{Gr}) appears to be a problem for us.
\end{remar}

\section*{Acknowledgements}
The second author (DT) thanks Haruzo Hida, Samir Siksek and Patrick Allen for many helpful discussions during the preparation of this article.
DT also acknowledges support of a Visiting Fellowship from ICMS Edinburgh,
which supported his work at the University of Warwick when this project was initiated.

\end{document}